\documentclass[12pt,a4paper]{article}

\usepackage[utf8]{inputenc}
\usepackage[english]{babel}
\usepackage{amsthm}
\usepackage{amsfonts}
\usepackage{amssymb}
\usepackage{amsmath}
\usepackage[affil-it]{authblk}
\usepackage{mathrsfs}
\usepackage{graphicx}
\usepackage{cite}

\newtheorem{Theorem}{Theorem}

\newtheorem{Lemma}{Lemma}

\newtheorem{Proposition}{Proposition}

\newenvironment{Corollary}{\par\noindent{\bf Corollary. }\it}{\par\bigskip}

\DeclareMathOperator{\re}{Re}
\DeclareMathOperator{\im}{Im}

\DeclareMathOperator{\loc}{loc}

\newcommand{\CC}{\mathbb{C}}
\newcommand{\RR}{\mathbb{R}}
\newcommand{\NN}{\mathbb{N}}

\newcommand{\const}{\mbox{const}}

\makeatletter
\renewcommand{\section}{\@startsection{section}{1}%
{0pt}{3.5ex plus 1ex minus .2ex}%
{2.3ex plus.2ex}{\normalfont\large\bfseries}}

\makeatother

\title{Completeness theorem for the system of eigenfunctions of the complex Schrödinger operator $\mathscr{L}_{c,\alpha}=-d^2/dx^2+cx^\alpha$}
\author{Sergey Tumanov
\thanks{sntumanov@yandex.ru}}
\affil{Moscow Center of Fundamental and Applied Mathematics at Lomonosov Moscow
State University, Russian Federation}

\begin{document}

\maketitle
\begin{abstract}
The completeness of the system of eigenfunctions of the complex Schrödinger
operator $\mathscr{L}_{c,\alpha}=-d^2/dx^2+cx^\alpha$ on the semi-axis with Dirichlet
boundary conditions is proved for an arbitrary $\alpha\in(0,2)$ and $|\arg c|<2\pi\alpha/(\alpha+2)+\Delta t(\alpha)$ with some $\Delta t(\alpha)>0$.
\end{abstract}
\frenchspacing

%------------------------------------------------------------
\section{Introduction}

We consider the operator
$$
\mathscr{L}_{c,\alpha}=-\frac{d^2}{dx^2}+cx^{\alpha}
$$
in $L_2(\RR_+)$ with Dirichlet boundary conditions with $c\in\CC$, $|\arg c|<\pi$, $\alpha>0$.

It's known, $\mathscr{L}_{c,\alpha}$ has a compact inverse, the spectrum is discrete, root subspaces are one-dimensional \cite{Tumanov23}.

For $0<|\arg c|<\pi$ it is not self-adjoint, moreover, it has bad spectral properties: the norm of the resolvent \cite{Davies, Krejc} and
the norms of spectral projectors \cite{Mityagin} grow exponentially.
Under these conditions, the operator cannot be similar to self-adjoint, its eigenfunctions do not form a Riesz basis in $L_2(\RR_+)$.
Nevertheless, the completeness of the system of eigenfunctions (S.E.) in $L_2(\RR_+)$ is an open problem, which is what our work is devoted to.

For $\alpha\ge2$ the problem of completeness of S.E. of $\mathscr{L}_{c,\alpha}$ is fully explored \cite{Davies,KrejcSiegl,Savchuk1}: the system is complete
for all $c\in\CC$: $|\arg c|<\pi$.

For $\alpha\in(0,2)$ completeness is proved in case $|\arg c|<t_0(\alpha)=2\pi\alpha/(\alpha+2)$ \cite{Savchuk1}. At the same time, the case
$t_0(\alpha)\le|\arg c|<\pi$ is almost not studied, since it is a much more difficult one.
The corresponding arguments are given in \cite{Tumanov23,Savchuk1}.

Perhaps the first result for the case $t_0(\alpha)\le|\arg c|<\pi$ was the study of Savchuk and Shkalikov \cite{Savchuk1} of the complex Airy operator ($\alpha=1$).
The authors proved the completeness of S.E. of $\mathscr{L}_{c,1}$ in case $|\arg c|<t_0(1)+\pi/6=5\pi/6$.

In our recent paper \cite{Tumanov23} a positive answer was given regarding the completeness of S.E. of $\mathscr{L}_{c,2/3}$
in case $|\arg c|<t_0(2/3)+\Delta t=\pi/2+\Delta t$, where $\Delta t>\pi/10$ is the only zero of some transcendental equation.

It turns out that these results generalize for all $\alpha\in(0,2)$: there exists $\Delta t=\Delta t(\alpha)>0$ continuously depending on $\alpha$,
such that the S.E. of $\mathscr{L}_{c,\alpha}$ is complete in case $|\arg c|<t_0(\alpha)+\Delta t(\alpha)$. This fact is the subject of our study.

The interest to the operators of the form $\mathscr{L}_{c,\alpha}$ has increased especially in recent decades due to the popularity of the ideas
of non-Hermitian quantum mechanics \cite{Bender1, Bender2}, as well as the general interest in spectral properties of
non-self-adjoint operators.

Let us note in this regard the problem posed by Almog in 2015 \cite{Almog}: is the S.E. of $\mathscr{L}_{i,\alpha}$ complete
for $0<\alpha\le2/3$? The answer is affirmative if $\theta_0(\alpha)=t_0(\alpha)+\Delta t(\alpha)>\pi/2$, in particular,
for $\alpha=2/3$. Taking into account the continuity of $\theta_0(\alpha)$,
there exists $\alpha_0<2/3$ such that the S.E. of $\mathscr{L}_{i,\alpha}$ is complete for all $\alpha>\alpha_0$. This hypothesis was
expressed by Savchuk and Shkalikov \cite{Savchuk1}. In our paper it is a simple Corollary to the main Theorem.
\bigskip

For complex $\zeta=|\zeta|e^{i\arg\zeta}$, $-\pi<\arg\zeta\le\pi$ and real $\beta$ we use the term $\zeta^\beta$ to denote the main branch:
$\zeta^\beta=|\zeta|^\beta e^{i\beta\arg\zeta}$.

Let $\theta\in[t_0(\alpha),\pi)\cap[t_0(\alpha),\pi\alpha)$, %]$
and
$$
\zeta_0(\theta)=e^{i(t_0(\alpha)-\theta)/\alpha},\quad Z_0(\theta)=(\sin t_0(\alpha)/\sin\theta)^{1/\alpha},
$$
we determine
\begin{gather*}
\rho(\theta)=\re\,\Bigl\{\int\limits_0^{Z_0(\theta)}\sqrt{e^{i\theta}\zeta^\alpha-e^{it_0(\alpha)}}\,d\zeta
-2\int\limits_0^{\zeta_0(\theta)}\sqrt{e^{i\theta}\zeta^\alpha-e^{it_0(\alpha)}}\,d\zeta\Bigr\}=\\
=\re\,\Bigl\{\int\limits_{\zeta_0(\theta)}^{Z_0(\theta)}\sqrt{e^{i\theta}\zeta^\alpha-e^{it_0(\alpha)}}\,d\zeta
-\int\limits_0^{\zeta_0(\theta)}\sqrt{e^{i\theta}\zeta^\alpha-e^{it_0(\alpha)}}\,d\zeta\Bigr\},
\end{gather*}
where the integration is performed over segments and the branch of the square root is chosen so that
$$
\re\int\limits_0^{Z_0(\theta)}\sqrt{e^{i\theta}\zeta^\alpha-e^{it_0(\alpha)}}\,d\zeta>0,\quad
\re\int\limits_0^{\zeta_0(\theta)}\sqrt{e^{i\theta}\zeta^\alpha-e^{it_0(\alpha)}}\,d\zeta>0.
$$

\begin{Theorem}
\label{th01}
Given $\alpha\in(0,2)$, the function $\rho(\theta)$ has the only zero $\theta_0(\alpha)$
within the interval:
$\theta_0(\alpha)\in (t_0(\alpha),\pi)\cap (t_0(\alpha),\pi\alpha)$; $\theta_0(\alpha)=t_0(\alpha)+\Delta t(\alpha)$,
$\Delta t(\alpha)>0$.

The function $\theta_0(\alpha)$ is continuous for $\alpha\in(0,2)$.

For $|\arg c|<\theta_0(\alpha)$ the S.E. of $\mathscr{L}_{c,\alpha}$ is complete in $L_2(\RR_+)$.
\end{Theorem}
\begin{Corollary}
There exists $\alpha_0<2/3$ such that the S.E. of $\mathscr{L}_{i,\alpha}$ is complete in $L_2(\RR_+)$
for all $\alpha>\alpha_0$.
\end{Corollary}

\begin{figure}
\begin{center}
\includegraphics[width=6cm,keepaspectratio]{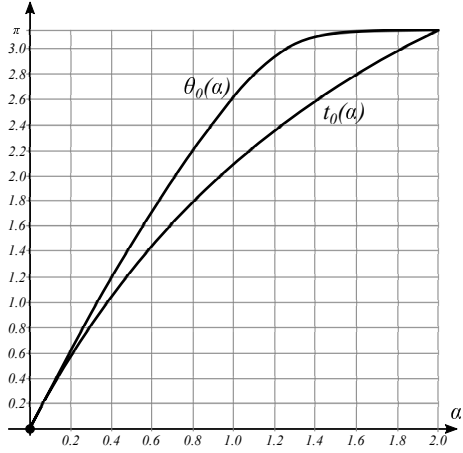}
\end{center}
\caption{Schematic chart $t_0(\alpha)$ and $\theta_0(\alpha)$.}
\label{theta0curves}
\end{figure}

The previously known completeness boundary $t_0(\alpha)$ and the one obtained in the present work $\theta_0(\alpha)$ are shown in the figure
\ref{theta0curves}.

\section{Auxiliary results}

The operator $\mathscr{L}_{c,\alpha}$ in $L_2(\RR_+)$ with Dirichlet boundary conditions is defined by differential expression
$$
l_{c,\alpha}(y)=-\frac{d^2y}{dx^2}+cx^\alpha y,\quad x\in[0,+\infty),%]
$$
and is considered on the domainn
$$
\mathfrak{D}(\mathscr{L}_{c,\alpha})=\bigl\{
y\in L_2(\RR_+)\,\bigl|\bigr.\, y\in W_{2,\loc}^2(\RR_+),\ l_{c,\alpha}(y)\in L_2(\RR_+),\ y(0)=0
\bigr\}.
$$

The following statements underlie the proof of the Theorem \ref{th01}. They correspond to Proposition A.1 and Lemmas 1 and 2 of \cite{Tumanov23}.

\begin{Lemma}
\label{lm01} Given $|\arg c|<\pi$, the operator $\mathscr{L}_{c,\alpha}$ is closed, has compact inverse.
The eigenvalues $\{\lambda_n\}$, $n\in\NN$ are simple (the root subspaces are one-dimensional), and have the form of
$\lambda_n=c^{2/(\alpha+2)}\tau_n$, $n\in\NN$, where $\tau_n>0$ don't depend on $c$.
\end{Lemma}

Consider the equation
\begin{equation}
y''=(cx^\alpha-\lambda)y.
\label{eqmaineqy}
\end{equation}

The existence of the so-called {\it Weyl solution} ---
the solution in $L_2(\RR_+)$ is known \cite[Ch II]{Titchmarsh} for any homogeneous Sturm--Liouville equation with real locally integrable potential.
This theory was generalized by Lidskii to the
case of some complex potentials \cite{Lidsky}. The following Lemma states for the equation \eqref{eqmaineqy}:
\begin{Lemma}[On the Weyl solution]
\label{lm02}
Given $0<|\arg c|<\pi$, there exists $\mathcal{Y}_0(x,\lambda)$ ($x\in\RR_+$, $\lambda\in\CC$) --- non-trivial solution to \eqref{eqmaineqy} with properties:
\begin{itemize}
\item the function $\mathcal{Y}_0(x,\lambda)$ is continuous of two variables $(x,\lambda)\in\RR_+\times\CC$;
\item for any $\lambda\in\CC$, $\mathcal{Y}_0(x,\lambda)\in L_2(\RR_+)$ as a function of $x$;
\item for any $x\ge0$, $\mathcal{Y}_0(x,\lambda)$ is an entire function of $\lambda$ of the order of growth of $\le(\alpha+2)/(2\alpha)$;
\item the zeros of $\mathcal{Y}_0(0,\lambda)$ coincide with the eigenvalues of $\mathscr{L}_{c,\alpha}$;
\item for an arbitrary $f\in L_2(\RR_+)$ the integral
\begin{equation}
\label{GfuncDefinition}
\mathcal{G}(\lambda)=\int\limits_0^{+\infty}\mathcal{Y}_0(x,\lambda)f(x)\,dx
\end{equation}
is an entire function of the order of growth of $\le(2+\alpha)/(2\alpha)$;
\item the function
$\mathcal{F}(\lambda)=\mathcal{G}(\lambda)/\mathcal{Y}_0(0,\lambda)$ is bounded in any closed sector for which
$\lambda\not\in\Lambda=\{\lambda\in\CC\,\Bigl|\Bigr.\,\arg\lambda\in[0,\arg c]\}$.
\end{itemize}
\end{Lemma}

Referring to Levin \cite{Levin}, the system of eigenfunctions $\{y_n\}$, $n\in\NN$ of the operator $\mathscr{L}_{c,\alpha}$
is generated by
$\mathcal{Y}_0(x,\lambda)$ and the sequence of eigenvalues $\{\lambda_n\}$: $y_n(x)=\mathcal{Y}_0(x,\lambda_n)$. It turns out that $\mathcal{Y}_0(x,\lambda)$ is
a closed kernel in terms of Levin: if for some $f\in L_2(\RR_+)$ holds $\mathcal{G}(\lambda)\equiv0$, then $f\equiv0$. Moreover,
\begin{Lemma}[On the closed kernel]
\label{lm03}
Given $0<|\arg c|<\pi$ and some $f\in L_2(\RR_+)$. If the function
\begin{equation}
\label{eqFFuncDef}
\mathcal{F}(\lambda)=\frac{1}{\mathcal{Y}_0(0,\lambda)}\int\limits_0^{+\infty}\mathcal{Y}_0(x,\lambda)f(x)\,dx
\end{equation}
is a constant, then $f\equiv0$.
\end{Lemma}

\section{The proof of the Theorem \ref{th01}}

In our proof we use some notation.

In estimates, $C$ denotes arbitrary positive constants, possibly different on both sides of the inequalities.

If for some $a>0$ we say that the statement is true for $a\ll1$ ($a\gg1$), we mean
the validity of this statement asymptotically --- for all $0<a<a_0$ ($a>a_0$) for some $a_0>0$.

Sometimes we use short notation: $\phi\equiv\arg\lambda$, $\theta\equiv\arg c$.

\bigskip

Let fix some $\alpha\in(0,2)$. Without loss of generality, $\im c>0$, $|c|=1$ --- with complex conjugation and scaling in $x$
one can cover any $c\in\CC$: $0<|\arg c|<\pi$.

Restrict ourselves to the case $\arg c\in [t_0(\alpha),\pi)\cap[t_0(\alpha),\pi\alpha)$. %]$

Let fix some $\delta>0$: $0<\delta<(\pi\min\{1,\alpha\}-\theta)/2$. Additional restrictions on $\delta>0$ will be introduced if necessary.

Let $\mathfrak{l}_\delta=\{\lambda\in\CC\,\bigl|\bigr.\,t_0(\alpha)-\delta<\arg\lambda<t_0(\alpha)\}$.
With our restrictions, $\arg c>\arg\lambda>0$ for $\lambda\in \mathfrak{l}_\delta$.
\bigskip

Let's walk through the main idea of the proof and its stages.

Assuming that for some $c$ the S.E. $\{y_n\}$, $n\in\NN$ is incomplete, one can find $f\not\equiv 0$, $f\in L_2(\RR_+)$, such that
$\mathcal{F}(\lambda)$ \eqref{eqFFuncDef} turns out to be an entire function. If we show that $\mathcal{F}(\lambda)$ is constant, then in view of the Lemma \ref{lm03}
we come to a contradiction with $f\not\equiv0$.

It is enough to prove that $\mathcal{F}(\lambda)$ grows not faster than a polynomial just in the small sector $\mathfrak{l}_\delta$. Indeed, on the one hand, the order
of growth of $\mathcal{F}(\lambda)$ is not higher than $\pi/t_0(\alpha)$, and on the other, ---
$\mathcal{F}(\lambda)$ is bounded in closed sectors outside
$\Lambda$ (Lemma \ref{lm02}). The central angles of the sectors complementing $\mathfrak{l}_\delta$ to $\Lambda$ are strictly
less than $t_0(\alpha)$, so referring to Phragm\'en-Lindel\"of Principle, we obtain that $\mathcal{F}(\lambda)$ is a polynomial.
Taking into account the boundedness of  $\mathcal{F}(\lambda)$ outside $\Lambda$, we obtain that $\mathcal{F}(\lambda)$ is a constant.

The condition $0<\theta<\theta_0(\alpha)$ turns out to be sufficient for such estimate of $\mathcal{F}(\lambda)$ in $\mathfrak{l}_\delta$.
\bigskip

We split the proof into several steps:
\begin{itemize}
\item {\bf Step 1}. The independent variable and parameter substitution in \eqref{eqmaineqy}.
\item {\bf Step 2}. Construction of the Weyl solution in terms of new variable and parameter, investigation of its uniform in $x\in\RR_+$ asymptotics for
$\mathfrak{l}_\delta\ni\lambda\to\infty$.
\item {\bf Step 3}. Proving that the condition $\rho(\arg c)<0$ is sufficient for the completeness of the S.E. of $\mathscr{L}_{c,\alpha}$.
\item {\bf Step 4}. Verification of the sufficient condition. The proof of the existence of a critical $\Delta t(\alpha)>0$
continuously depending on $\alpha$ such that the condition $\rho(\theta)<0$ is satisfied
for all $\theta$: $t_0(\alpha)\le\theta<t_0(\alpha)+\Delta t(\alpha)$. This completes the proof of the Theorem \ref{th01}.
\end{itemize}

\subsection{Step 1. The independent variable and parameter substitution}

Let $\lambda\in\mathfrak{l}_\delta$. We fix some $\arg\lambda$, at the same time $|\lambda|>0$ can be arbitrarily large.

For any solution $y(x)=y(x,\lambda)$
to the equation \eqref{eqmaineqy}, let
$$
w(t)=y(|\lambda|^{1/\alpha}t),\quad k=|\lambda|^{1/2+1/\alpha}.
$$
The equation is converted to the form:
\begin{equation}
w''=k^2(ct^\alpha-\mu)w,
\label{eqmaineqw}
\end{equation}
with $\mu=e^{i\phi}$, $c=e^{i\theta}$, $0<\phi<\theta$.

\begin{Proposition}
\label{prop01}
The following inequalities hold:
\begin{equation}
0<\arg c/\mu<\frac{\pi}{3},\quad
0<\arg (c/\mu)^{1/\alpha}<\frac{\pi}{3}.
\label{eqlamcdiv}
\end{equation}
\end{Proposition}
{\noindent\bf Proof.} Only the upper bounds are of interest.

For $0<\alpha<1$ holds $\delta+\arg c<\pi\alpha$. Thus
\begin{equation*}
\begin{split}
\arg c/\mu &<\pi\alpha-t_0(\alpha)=\frac{\pi\alpha^2}{\alpha+2}<\pi/3,\\
\arg (c/\mu)^{1/\alpha} &<\frac{\pi\alpha}{\alpha+2}<\pi/3,
\end{split}
\end{equation*}

For $1\le\alpha<2$ holds $\delta+\arg c<\pi$. Thus
\begin{equation*}
\begin{split}
\arg c/\mu &<\pi-t_0(\alpha)=\frac{\pi(2-\alpha)}{\alpha+2}\le\pi/3,\\
\arg (c/\mu)^{1/\alpha} &<\frac{\pi(2-\alpha)}{\alpha(\alpha+2)}\le\pi/3.
\end{split}
\end{equation*}
In each case, the estimates reach their maximum at $\alpha=1$.\qquad$\Box$

\subsection{Step 2. Construction of the Weyl solution}
Further, let $\mu\in\mathfrak{l}_\delta$, $\mu=e^{i\phi}$.

Determine a multivalued analytic function in $\CC_r=\{z\in\CC\,|\,\re z>0\}$:
$$
S(z)=\int\limits_0^z\sqrt{c\zeta^\alpha-\mu}\,d\zeta,
$$
integrating over rectifiable paths in the right half-plane that
do not pass through the only singularity $\zeta_0=(\mu/c)^{1/\alpha}$ in the IVth quarter \eqref{eqlamcdiv}.

\begin{figure}
\begin{center}
\includegraphics[width=6cm,keepaspectratio]{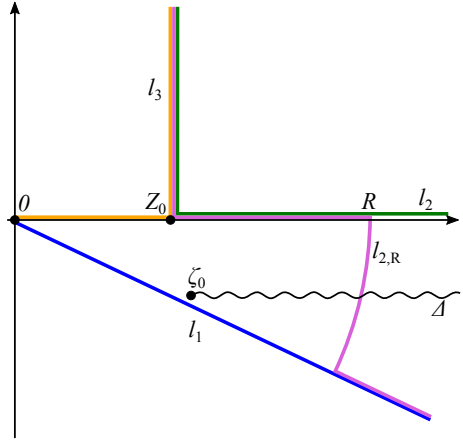}
\end{center}
\caption{The paths of monotonicity of $\re S(z)$. The wavy line corresponds to the cut.}
\label{figljcurves}
\end{figure}

Determine the horizontal cut $\Delta=\{z=\zeta_0+t,\ t\ge0\}$ (see fig.\ref{figljcurves}).

For $z\in\CC_r$ we define $q(z)=cz^\alpha-\mu$.

We also define the main branches of $q^{1/2}(z)$ and $S(z)$ in the simply connected domain $\CC_r\setminus\Delta$ by the conditions:
$S(0)=0$ (by continuity); $\re q^{1/2}(z)>0$ and $\re S(z)>0$ for $z>0$, $z\ll1$.
In what follows, $\sqrt{cz^\alpha-\mu}$ and $S(z)$ will denote
the main branches of $q^{1/2}(z)$ and $S(z)$.
Cases of the continuation of these branches through the cut $\Delta$ will be specially noted.

We also fix in $\CC_r\setminus\Delta$ any branch of $q^{1/4}(z)$, which we call the main one.

Denote $Z_0=(\im\mu/\im c)^{1/\alpha}>0$.

\begin{Proposition}
\label{propSonLs}
There exist $\varepsilon>0$, $R_0>0$, the value $\re S(z)$ is monotone along each of the following paths for $R>R_0$:
\begin{gather*}
l_1=\{\zeta_0e^{-i\varepsilon}t,\ 0\le t<+\infty\},\\
l_2=\{Z_0+it,\ t\ge0\}\cup[Z_0,+\infty),\\%]
l_{2,R}=\{Z_0+it,\ t\ge0\}\cup[Z_0,R]\cup\Gamma_R\cup\{z\in l_1,\ |z|\ge R\},\\
l_3=[0,Z_0]\cup\{Z_0+it,\ t\ge0\},
\end{gather*}
where $\Gamma_R=\{Re^{it},\ t\in[\arg\zeta_0-\varepsilon,0]\}$ (see fig.\ref{figljcurves}).
The analytic continuation of the main branch of $S(z)$ is considered along $l_{2,R}$, since $l_{2,R}$ is crossing $\Delta$.

When $z$ moves from the origin to an infinite point along each of the paths $l_1$ and $l_3$, the value of $\re S(z)$ increases monotonically from $0$ to $+\infty$.

When $z$ moves from $+\infty$ to $Z_0+i\infty$ along $l_2$, the value of $\re S(z)$ increases monotonically from $-\infty$ to $+\infty$.

When $z$ moves from $\zeta_0e^{-i\varepsilon}\infty$ to $Z_0+i\infty$ along $l_{2,R}$, the value of $\re S(z)$ decreases monotonically from $+\infty$ to $-\infty$.
The continuation of $S(z)$ through $\Delta$ is considered here.

Along the real axis the value $\re S(z)$: is increasing on $[0,Z_0]$, has a single maximum at $z=Z_0$, is decreasing to $-\infty$ on $(Z_0,+\infty)$
with asymptotics:
\begin{equation}
\label{eqSonRplus}
S(z)\sim -\ \frac{2\re c^{1/2}}{\alpha+2}z^{\alpha/2+1},\mbox{ as }z\to+\infty.
\end{equation}

For all $z>0$ the value $\im q^{1/2}(z)<0$ for the main branch of $q^{1/2}$.

The inequality holds: $\re S(\zeta_0)>0$.
\end{Proposition}
{\noindent\bf Proof.} We use the asymptotic formula following from the integral representation of $S(z)$ in neighborhoods of the infinite parts of the curves $l_j$:
\begin{equation}
S(z)\sim\int\limits_0^z\sqrt{c\zeta^\alpha}\,d\zeta\sim\pm c^{1/2}\frac{z^{\alpha/2+1}}{\alpha/2+1},\quad z\to\infty,
\label{eqSinInfAsympt}
\end{equation}
the choice of the sign will be further specified for each curve $l_j$.

It also follows from the integral representation of $S(z)$:
\begin{equation}
S(z)\sim\int\limits_0^z\sqrt{-\mu}\,d\zeta\sim -i\mu^{1/2}z\mbox{, as }z\to0,\ \re z>0,
\label{eqSZeroAssymp}
\end{equation}
the choice of the sign is determined by the condition $\re S(z)>0$ for $z>0$, $z\ll1$.

It follows from \eqref{eqSZeroAssymp} that $\im S(z)<0$ for $z\ll1$.

Further we explore the extrema of $\re S(z)$ along the real axis for $z=t\ge0$:
$$
\frac{d}{dt}\re S(t)=\re\sqrt{ct^\alpha-\mu}=0\ \Leftrightarrow\
\sqrt{ct^\alpha-\mu}=i\beta_1
$$
for some $t\ge0$, $\beta_1\in\RR\setminus\{0\}$, that is equivalent to $ct^\alpha-\mu=-\beta_1^2<0$. Considering the imaginary part,
$t=Z_0=(\im\mu/\im c)^{1/\alpha}$.
For the real part,
$$
\re cZ_0^\alpha-\re\mu=\frac{\im\{\mu \overline c\}}{\im c}=\frac{|c|^2}{\im c}\im\frac{\mu}{c}<0,
$$
taking \eqref{eqlamcdiv} into account. Thus there is a unique extremum of $\re S(z)$ at $Z_0>0$ --- the global maximum for $z\ge0$.

Exploring the extrema of $\im S(z)$ for $z=t\ge0$:
$$
\frac{d}{dt}\im S(t)=0\ \Leftrightarrow\
\sqrt{ct^\alpha-\mu}=\beta_2\in\RR\setminus\{0\},
$$
that is equivalent to $ct^\alpha-\mu=\beta_2^2>0$. Considering the imaginary part, again we get $t=Z_0$, but $cZ_0^\alpha-\mu=-\beta_1^2<0$, thus,
$\im S(z)$ has no extrema for $z\ge0$.

To summarize: for $z\ge0$ the function $\re S(z)$ increases, reaches its maximum at $z=Z_0$, then decreases. The function $\im S(z)$ decreases for all $z\ge0$.
In particular for the main branch of $q^{1/2}(z)$ we obtain $\im q^{1/2}(z)=\im S'(z)<0$ for $z> 0$.

Further we explore the extrema of $\re S(z)$ along the vertical rays $\{a+it,\ t\ge0\}$ for $a>0$:
$$
\frac{d}{dt}\re S(a+it)=0\ \Leftrightarrow\
c(a+it)^\alpha-\mu=\beta_3^2>0\ \Leftrightarrow\
(a+it)^\alpha=\frac{\mu}{c}+\frac{\beta_3^2}{c}.
$$

Both terms on the right-hand side lie in the lower half-plane \eqref{eqlamcdiv}. At the same time, since $\alpha\in(0,2)$,
the value $(a+it)^\alpha$
lies in the upper half-plane, or on the real axis. In other words, there are no extrema, $\re S(z)$ is strictly monotone. In particular, along the
vertical ray $l_3\cap l_2$.

Let's move along $l_3$ from the origin and observe the curve in the image of $S(z)$. The value $\re S(z)$ is strictly increasing, and the value $\im S(z)$
is strictly decreasing on the segment $[0,Z_0]$. The tangent to the image $S([0,Z_0])$ at $S(Z_0)$ is strictly vertical and directed downward.
Since $S(z)$ is univalent in $Z_0$, the angles are preserved: when moving  along $l_3$ from $Z_0$ to $Z_0+i\infty$, the value $\re S(z)$ will increase.

The same reasoning shows that $\re S(z)$ is monotone along $l_2$.

We take $\varepsilon>0$ so that $l_1$ lies in the IVth quarter. We reduce $\varepsilon>0$ so that $S(z)$ lies in the IVth quarter for arbitrarily small points of
$l_1$. This is possible, as $\delta<(\pi\alpha-\theta)/2$. Indeed,
\begin{gather}
\notag
\lim_{t\to+0}\arg S(\zeta_0e^{-i\varepsilon}t)=\arg \bigl\{-i\zeta_0e^{-i\varepsilon}\mu^{1/2}\bigr\}=\\
\notag
=-\frac{\pi}{2}+\frac{t_0(\alpha)}{2}-\frac{\delta_1}{2}+
\frac{t_0(\alpha)}{\alpha}-\frac{\delta_1}{\alpha}-\frac{\theta}{\alpha}-\varepsilon
=\frac{\pi}{2}-\frac{\theta}{\alpha}-\delta_1\frac{\alpha+2}{2\alpha}-\varepsilon>\\
\label{eqargSonl1estim}
>\frac{\pi}{2}-\frac{\theta}{\alpha}-\delta\frac{\alpha+2}{2\alpha}-\varepsilon
\end{gather}
where $\delta_1=t_0(\alpha)-\phi\in(0,\delta)$. Since $\alpha\in(0,2)$,
$$
\frac{\pi}{2}-\frac{\theta}{\alpha}-\delta\frac{\alpha+2}{2\alpha}>
\frac{\pi}{2}-\frac{\theta}{\alpha}-\frac{1}{2}(\pi\alpha-\theta)\frac{\alpha+2}{2\alpha}>
\frac{\pi}{2}-\frac{\theta}{\alpha}-\frac{2}{\alpha+2}(\pi\alpha-\theta)\frac{\alpha+2}{2\alpha}=-\frac{\pi}{2},
$$
therefore, we reduce $\varepsilon>0$ to make the final estimate \eqref{eqargSonl1estim} greater than $-\pi/2$.
On the other hand, as $\pi/2<t_0(\alpha)/\alpha\le\theta/\alpha$, we have
$$
\lim_{t\to+0}\arg S(\zeta_0e^{-i\varepsilon}t)=\frac{\pi}{2}-\frac{\theta}{\alpha}-\delta_1\frac{\alpha+2}{2\alpha}-\varepsilon<0.
$$

Further we explore the extrema of $\re S(z)$ along $l_1$:
\begin{equation}
\frac{d}{dt}\re S(\zeta_0e^{-i\varepsilon}t)=0\ \Leftrightarrow\
\re\bigl\{\zeta_0e^{-i\varepsilon}\mu^{1/2}\sqrt{e^{-i\varepsilon\alpha}t^\alpha-1}\bigr\}=0.
\label{eqExtrOnl1}
\end{equation}

The value $\zeta_0e^{-i\varepsilon}\mu^{1/2}$ lies in the Ist quarter (see the first equality in \eqref{eqargSonl1estim});
for $t>0$, the expression under the root in \eqref{eqExtrOnl1} lies is in the lower half-plane, that is, the root itself is in the IVth quarter (up to a sign). The whole
expression in curly braces in \eqref{eqExtrOnl1} lies in the right half-plane (up to a sign), i.e. equality \eqref{eqExtrOnl1} is impossible,
and $\re S(z)$ has no extrema along $l_1$.

With $\varepsilon=0$, the same reasoning proves the monotonic growth of $\re S(\zeta_0t)$, along $0\le t<1$. Whence $\re S(\zeta_0)>0$.

Now we show that  $\re S(z)$ is monotonic along $\Gamma_R=\{Re^{it},\ t\in[\arg\zeta_0-\varepsilon,0]\}$:
\begin{equation}
\frac{d}{dt}\re S(Re^{it})=0\ \Leftrightarrow\
e^{2it}(cR^\alpha e^{it\alpha}-\mu)=\beta_4^2>0.
\label{eqExtrOnGammaR}
\end{equation}

For $R\gg1$, the argument of the left part of the last equality in \eqref{eqExtrOnGammaR} is determined by the argument of $ce^{it(\alpha+2)}$, but
\begin{gather}
\notag
\arg \{ce^{it(\alpha+2)}\}=\theta+t(\alpha+2)\ge\theta+\bigl(\frac{\phi-\theta}{\alpha}-\varepsilon\bigr)(\alpha+2)=\\
\label{eqOnGammaEvl}
=2\frac{\pi\alpha-\theta}{\alpha}-\delta_1\frac{\alpha+2}{\alpha}-\varepsilon(\alpha+2),
\end{gather}
where $\delta_1=t_0(\alpha)-\phi\in(0,\delta)$. As $\delta<(\pi\alpha-\theta)/2$, we have:
\begin{gather*}
2\frac{\pi\alpha-\theta}{\alpha}-\delta_1\frac{\alpha+2}{\alpha}>2\frac{\pi\alpha-\theta}{\alpha}-\delta\frac{\alpha+2}{\alpha}=
\frac{\alpha+2}{\alpha}\Bigl(
\frac{2}{\alpha+2}(\pi\alpha-\theta)-\delta\Bigr)>\\
>\frac{\alpha+2}{\alpha}\Bigl(
\frac{1}{2}(\pi\alpha-\theta)-\delta
\Bigr)>0,
\end{gather*}
reducing $\varepsilon>0$ if needed, we achieve $\arg \{ce^{it(\alpha+2)}\}>0$ \eqref{eqOnGammaEvl}.

On the other hand, $\arg \{ce^{it(\alpha+2)}\}\le\arg c<\pi$ for $t\le0$, i.e. $ce^{it(\alpha+2)}$ lies in the upper half-plane, therefore, for sufficient
large $R>R_0>0$ the left part of the last equality in \eqref{eqExtrOnGammaR} also lies
strictly in the upper half-plane, i.e. cannot be real.

For $z=Re^{it}\in\Gamma_R$ the value $ce^{it(\alpha+2)}$ lies in the upper half-plane, $c^{1/2}e^{it(\alpha/2+1)}$ --- in the Ist quarter, in particular
for $t=\arg\zeta_0-\varepsilon$. Thus, for $z\in l_1$, $z=|z|e^{it}$, $|z|\to\infty$:
$$
S(z)\sim c^{1/2}e^{it(\alpha/2+1)}\frac{|z|^{\alpha/2+1}}{\alpha/2+1},
$$
the value $S(z)$ has to lie in the right half-plane due to the monotonic growth of $\re S(z)$ along $l_1$, this determines the choice of the sign in the last formula,
and $\re S(z)\to+\infty$ as $l_1\ni z\to\infty$.

Let $z\in\Gamma_R$, consider the analytic continuation of $S(z)$ through the cut $\Delta$ from below, \eqref{eqSinInfAsympt} takes the form:
\begin{equation}
S(z)\sim c^{1/2}\frac{z^{\alpha/2+1}}{\alpha/2+1},\mbox{ uniformly in } z\in\Gamma_R,\ R\to+\infty.
\label{eqSinMAsympt}
\end{equation}
For $R\gg1$, $z\in\Gamma_R$ the value $\re S(z)$ lies in the Ist quarter. Taking into account the monotonicity of $\re S(z)$ along $\Gamma_R$, it follows from \eqref{eqSinMAsympt}
that $\re S(z)$ decreases with $z$ moving along $\Gamma_R$ counterclock\-wise for $R\gg1$.
In view of the previous considerations along the rest parts of $l_{2,R}$, this justifies the monotonicity of $\re S(z)$
along the entire $l_{2,R}$ for $R\gg1$, without loss of generality, for $R>R_0>0$.

The formula \eqref{eqSonRplus} follows from \eqref{eqSinInfAsympt} up to a sign. We conclude that $\re S(z)$ is unbounded for $z\to+\infty$.
Since $\re S(z)$ decreases on $z\in[Z_0,+\infty)$, %]$
we obtain $\re S(z)\to-\infty$, and the sign ``minus'' is selected in \eqref{eqSonRplus}.

Let $z\in l_3$, $z\to\infty$. Then $z=|z|e^{i\pi/2+o(1)}\sim|z|e^{i\pi/2}$. Using the asymptotics of the integral $S(z)$, we obtain:
$$
S(z)\sim -c^{1/2}e^{i\pi(\alpha/2+1)/2}\frac{|z|^{\alpha/2+1}}{\alpha/2+1},
$$
the term $e^{i\pi(\alpha/2+1)/2}$ lies in the IInd quarter, as $\alpha\in(0,2)$; the term $c^{1/2}$ lies in the Ist quarter.
The product $c^{1/2}e^{i\pi(\alpha/2+1)/2}$ lies in the left half-plane, thus $\re S(z)$ is unbounded, $\re S(z)\to+\infty$,
and the sign ``minus'' is used in the last asymptotic formula.\qquad$\Box$
\bigskip

The next Proposition is the main one in the current section.
\begin{Proposition}
\label{prop3}
For $k\gg1$ and $t\ge0$ there is the solution $\mathcal{W}_0(t,k)$ to \eqref{eqmaineqw}, for which the uniform asymptotic formulas are valid:
\begin{equation}
\label{eqWonRplus}
\begin{split}
\mathcal{W}_0(t,k)&=(c t^\alpha-\mu)^{-1/4}\bigl\{-ie^{k(S(t)-2S(\zeta_0))}[1]+e^{-kS(t)}[1]\bigr\},\quad t\in[0,Z_0],\\
\mathcal{W}_0(t,k)&=(c t^\alpha-\mu)^{-1/4}\bigl\{-ie^{k(S(t)-2S(\zeta_0))}[1]\bigr\},\quad t>Z_0,
\end{split}
\end{equation}
where $[1]=1+O(k^{-1})$ is uniform in $t$ on the corresponding interval.

For $k\gg1$, $\mathcal{W}_0(t,k)\in L_2(\RR_+)$ as a function of $t$.
\end{Proposition}
{\noindent\bf Proof.} Further, let $z\in\CC_r\cup\{0\}$.

The function $\re S(z)$ increases along $l_1$. For $k\gg1$, we construct the {\it subordinate} solution to \eqref{eqmaineqw}
at the infinite point of $l_1$, according to \cite[Ch.6 \S12]{Olver}, of the form
\begin{equation}
\mathcal{W}_0(z,k)=(c z^\alpha-\mu)^{-1/4}e^{-kS(z)}(1+\epsilon(z,k)),
\label{eqwasympl1}
\end{equation}
where
$$
|\epsilon(z,k)|\le\exp\Bigl\{\frac{\mathscr{V}_{z,\infty}(F)}{2k}\Bigr\}-1,\ F=\int\frac{1}{q^{1/4}}\frac{d^2}{dz^2}\Bigl(\frac{1}{q^{1/4}}\Bigr)\,dz.
$$
The symbol $\mathscr{V}_{z,\infty}(F)$ denotes to the variation of the {\it error-control function} $F(z)$ along the unbounded part of $l_1$,
starting from $z\in l_1$.

We estimate:
$$
\mathscr{V}_{z,\infty}(F)=\int\limits_{[z,\infty)\subset l_1}
\Bigl|\frac{1}{q^{1/4}}\frac{d^2}{dz^2}\Bigl(\frac{1}{q^{1/4}}\Bigr)\Bigr|\,|d\zeta|\le C
\int\limits_{[z,\infty)\subset l_1}\frac{|d\zeta|}{1+|\zeta|^{2+\alpha/2}}\le\frac{C}{1+|z|^{1+\alpha/2}},
$$
thus $\epsilon(z,k)=O(k^{-1})$ uniformly in $z\in l_1$.

According to \cite[Ch.5 \S3]{Olver}, the solution $\mathcal{W}_0(z,k)$ can be continued analytically to $z\in\CC_r$.

Consider $l_{2,R}$ for $R>R_0>0$ (Proposition \ref{propSonLs}) and construct the solution to \eqref{eqmaineqw} along $l_{2,R}$ subordinate
at the infinite point of the ray $l_1\cap l_{2,R}$.
By construction, it coincides with $\mathcal{W}_0(z,k)$ along $l_1\cap l_{2,R}$, thus it coincides with $\mathcal{W}_0(z,k)$
in $\CC_r$ due to the uniqueness of the analytic continuation.

The formula \eqref{eqwasympl1} remains along the entire $l_{2,R}$ with the remark, that when passing through the cut $\Delta$ from below, one should
consider the analytic continuation of the main branches $q^{1/4}(z)$ and $S(z)$.

The new branches of $q^{1/4}(z)$ and $S(z)$ obtained as a result of this continuation
will be also considered in the domain $\CC_r\setminus\Delta$ (on the second sheet). Eventually, for $z\in\CC_r\setminus\Delta$,
the new branch of $q^{1/4}(z)$ takes the form of $iq^{1/4}(z)$.
The new branch of $S(z)$ (denoted as $\tilde S(z)$) takes the form of:
\begin{equation}
\tilde S(z)=2S(\zeta_0)-S(z).
\label{eqtildeS}
\end{equation}

In terms of new branches, for $z\in (l_3\cap l_{2,R})\cup[Z_0,R]$, $k\gg1$:
\begin{gather}
\label{eqwasympl2R}
\mathcal{W}_0(z,k)=-i(c z^\alpha-\mu)^{-1/4}e^{-k\tilde S(z)}(1+\epsilon(z,k)),\\
\notag
|\epsilon(z,k)|\le\exp\Bigl\{\frac{\mathscr{V}_{z,\infty}(F)}{2k}\Bigr\}-1,
\end{gather}
where the variation is considered along the unbounded part of $l_{2,R}$ with the end at the infinite point of the ray $l_1\cap l_{2,R}$.
Obviously:
\begin{gather*}
\mathscr{V}_{z,\infty}(F)<C,\mbox{ for }z\in l_{2,R},\\
\mathscr{V}_{z,\infty}(F)<C\Bigl(\int\limits_z^R+\int\limits_{\Gamma_R}+\int\limits_{[z,\infty)\subset l_1}\Bigr)\frac{|d\zeta|}{|\zeta|^{2+\alpha/2}}\le
\frac{C}{z^{1+\alpha/2}}+\frac{C}{R^{1+\alpha/2}},\mbox{ for }z\in [Z_0,R],
\end{gather*}
where the constant $C>0$ does not depend on $R>R_0$.

Since the choice of $R>R_0$ is arbitrary, the formula \eqref{eqwasympl2R} defines a representation for $\mathcal{W}_0(z,k)$ along $l_2$
for $k\gg1$ with the uniform in $z\in l_2$ estimate $\epsilon(z,k)=O(k^{-1})$. This proves the second formula in \eqref{eqWonRplus}.

We turn to the path $l_3$, along which we construct two solutions to \eqref{eqmaineqw}:
$u(z,k)$ --- subordinate at $z=0$ and $v(z,k)$ --- subordinate at the infinite point of $l_3$.
Like $\mathcal{W}_0(z,k)$, both $u(z,k)$ and $v(z,k)$ can be analytically continued to $\CC_r$.
The following asymptotics are valid for $k\gg1$ uniformly in $z\in l_3$:
\begin{equation}
\begin{split}
u(z,k)&=(c z^\alpha-\mu)^{-1/4}e^{kS(z)}(1+O(k^{-1})),\\
v(z,k)&=(c z^\alpha-\mu)^{-1/4}e^{-kS(z)}(1+O(k^{-1})).
\end{split}
\label{equvasympl3}
\end{equation}
Moreover, \eqref{equvasympl3} guarantees the linear independence of $u(z,k)$ and $v(z,k)$ for $k\gg1$.

The constructed functions are solutions to a homogeneous second-order equation. For each $k\gg1$ there exist two constants $A(k)$ and $B(k)$:
\begin{equation}
\mathcal{W}_0(z,k)=A(k)u(z,k)+B(k)v(z,k).
\label{eqwviauv}
\end{equation}
To find $A(k)$ and $B(k)$, we use \eqref{eqwasympl1}, \eqref{eqwasympl2R},
\eqref{equvasympl3} at $z=0$ and at some arbitrary large in absolute value point $z_*\in\{Z_0+it,\ t\ge0\}\subset l_3$:
$$
\left\{
\begin{aligned}
A(k)[1]+B(k)[1]&=1\\
A(k)e^{kS(z_*)}[1]+B(k)e^{-kS(z_*)}[1]&=-ie^{-k\tilde S(z_*)}
\end{aligned}
\right.,
$$
where $[1]=1+O(k^{-1})$ for $k\gg1$ uniformly in $z\in\{Z_0+it,\ t\ge0\}$.

Taking into account \eqref{eqtildeS}, the limit $\re S(z)\to+\infty$ for $l_3\ni z\to\infty$, the inequality $\re S(\zeta_0)>0$ (Proposition \ref{propSonLs}),
finally:
$$
A(k)=-ie^{-2kS(\zeta_0)}[1],\quad
B(k)=[1].
$$

Combining this result with the formula \eqref{eqwviauv} and the asymptotics \eqref{equvasympl3}, we get the first formula in \eqref{eqWonRplus}.

For $k\gg1$, taking into account \eqref{eqSonRplus} and \eqref{eqWonRplus}, we obtain that $\mathcal{W}_0(t,k)$ is bounded and exponentially decaying as $t\to+\infty$,
thus $\mathcal{W}_0(t,k)\in L_2(\RR_+)$ as a function of $t$.\qquad$\Box$

\subsection{Step 3. The sufficient condition for the completeness}

Note that the function $S(z)$, the critical points $Z_0=(\im\mu/\im c)^{1/\alpha}>0$ and $\zeta_0=(\mu/c)^{1/\alpha}$
depend on $\theta=\arg c\in[t_0(\alpha),\pi)\cap[t_0(\alpha),\pi\alpha)$, %]$.
and on $\phi=\arg\mu\in(t_0(\alpha)-\delta,t_0(\alpha))$.

We determine:
$$
\rho(\theta,\phi)=\re\bigl\{S(Z_0)-2S(\zeta_0)\bigr\}.
$$

For each $\theta\in[t_0(\alpha),\pi)\cap[t_0(\alpha),\pi\alpha)$ %]$
the function $\rho(\theta,\phi)$ can be continuously extender to the point $\phi=t_0(\alpha)$. Determine
$$
\rho(\theta)=\rho(\theta,t_0(\alpha))=\re\bigl\{S(Z_0)-2S(\zeta_0)\bigr\}\Bigl|_{\phi=t_0(\alpha)}\Bigr..
$$

\begin{Proposition}
\label{propSuff}
The condition $\rho(\theta)<0$, $\theta=\arg c\in[t_0(\alpha),\pi)\cap[t_0(\alpha),\pi\alpha)$ %]$
is sufficient for the completeness of S.E. of the operator $\mathscr{L}_{c,\alpha}$.
\end{Proposition}
{\noindent\bf Proof.} Let $\rho(\theta)<0$. We take $\delta>0$, $\delta<(\pi\min\{1,\alpha\}-\theta)/2$ to satisfy $\rho(\theta,\phi)<0$ for all
$\mu\in\mathfrak{l}_\delta$.

Recently with the substitution of the variable and the parameter, we constructed the solution $\mathcal{W}_0(t,k)\in L_2(\RR_+)$ to
\eqref{eqmaineqw} for $k\gg1$. The function $\mathcal{Y}_0(|\lambda|^{1/\alpha}t,\lambda)\in L_2(\RR_+)$ is also the solution to \eqref{eqmaineqw}. Hereinafter
$k=|\lambda|^{1/2+1/\alpha}$.

The eigenvalues of $\mathscr{L}_{c,\alpha}$ lie on the ray $\arg\lambda=2\theta/(\alpha+2)$ (Lemma \ref{lm01}).
As $\delta<(\pi\alpha-\theta)/2$, we have: $2\theta/(\alpha+2)<t_0(\alpha)-\delta$, thus, there are no eigenvalues inside $\mathfrak{l}_\delta$.
Therefore for $\lambda\in\mathfrak{l}_\delta$, the solution $y(x)$ to \eqref{eqmaineqy} with the initial conditions $y(0)=0$, $y'(0)=1$
is not in $L_2(\RR_+)$; the subspace of solutions in
$L_2(\RR_+)$ is one-dimensional. The same is true for the equation \eqref{eqmaineqw}. Thus for $k\gg1$, there exists
the function $A(k)$ such that
\begin{equation}
\label{eqY0viaW}
\mathcal{Y}_0(|\lambda|^{1/\alpha}t,\lambda)=A(k)\mathcal{W}_0(t,k).
\end{equation}

Let us assume the contrary, that the S.E. of $\mathscr{L}_{c,\alpha}$ is incomplete. Then there exists $f\in L_2(\RR_+)$, $f\not\equiv0$, and
for all eigenvalues $\{\lambda_n\}$ of the operator $\mathscr{L}_{c,\alpha}$,
$\mathcal{G}(\lambda_n)=0$ \eqref{GfuncDefinition}. In this case $\mathcal{F}(\lambda)$ \eqref{eqFFuncDef} is an entire function.
We estimate it for $|\lambda|\gg1$ on an arbitrary ray $\mathfrak{l}\subset\mathfrak{l}_\delta$ with fixed $\phi=\arg\lambda$.
\begin{gather*}
|\mathcal{F}(\lambda)|^2\le C\frac{1}{|\mathcal{Y}_0(0,\lambda)|^2}\int\limits_0^{+\infty}|\mathcal{Y}_0(x,\lambda)|^2\,dx=
C\frac{1}{|\mathcal{Y}_0(0,\lambda)|^2}|\lambda|^{1/\alpha}\int\limits_0^{+\infty}|\mathcal{Y}_0(|\lambda|^{1/\alpha}t,\lambda)|^2\,dt=\\
=C\frac{1}{|\mathcal{Y}_0(0,\lambda)|^2}|\lambda|^{1/\alpha}\int\limits_0^{+\infty}|A(k)\mathcal{W}_0(t,k)|^2\,dt=
C\frac{1}{|\mathcal{W}_0(0,k)|^2}|\lambda|^{1/\alpha}\int\limits_0^{+\infty}|\mathcal{W}_0(t,k)|^2\,dt\le\\
\le C|\lambda|^{1/\alpha}\Bigl\{\int\limits_0^{Z_0}(e^{2k\re(S(t)-2S(\zeta_0))}+e^{-2k\re S(t)})\,dt+\int\limits_{Z_0}^{+\infty}e^{2k\re(S(t)-2S(\zeta_0))}\,dt\Bigr\}\le\\
\le C|\lambda|^{1/\alpha}\bigl\{e^{2k\re(S(Z_0)-2S(\zeta_0))}+1\bigr\}=C|\lambda|^{1/\alpha}\bigl\{e^{2k\rho(\theta,\phi)}+1\bigr\}<C|\lambda|^{1/\alpha}.
\end{gather*}

The first estimation follows from the Cauchy--Bunyakovsky inequality and the condition $f\in L_2(\RR_+)$; then we make the substitution
$x=|\lambda|^\alpha t$; then we use
\eqref{eqY0viaW} and boundedness of $|ct^\alpha-\mu|^{-1/4}$ along $\mathfrak{l}$. We take into account that $\re S(t)\ge 0$ for $t\in[0,Z_0]$. The final estimation
follows from Laplace method. Everywhere $C$ does not depend on $k$.

By construction $0<t_0(\alpha)-\delta<\phi<t_0(\alpha)\le\theta$. As $\delta<\pi\alpha-\theta$, it follows: $\delta<4\pi\alpha/(\alpha+2)-\theta$,
thus $\theta-t_0(\alpha)+\delta<t_0(\alpha)$.

In view of these estimates, we can form two adjacent sectors $\Lambda_1$, $\Lambda_2$ with a common boundary passing along the ray $\mathfrak{l}$, so that
1) the central angle of each sector is
less than $t_0(\alpha)$, 2) the second boundary of each sector lies outside $\Lambda=\{\arg\lambda\in[0,\theta]\}$.

The order of growth of $\mathcal{F}(\lambda)$ is at most $\pi/t_0(\alpha)$ (Lemma \ref{lm02}).
The function $\mathcal{F}(\lambda)$ is bounded on the rays $\arg\lambda\not\in\Lambda$
and grows no faster than a polynomial on the ray $\mathfrak{l}$. It follows from the Phragm\'en-Lindel\"of Principle that in each of the sectors
$\Lambda_j$ ($j=1,2$) the function $\mathcal{F}(\lambda)$ grows no faster than a polynomial, and therefore in the whole plane $\CC=\Lambda\cup\Lambda_1\cup\Lambda_2$.
In this way,
the entire function $\mathcal{F}(\lambda)$ is a polynomial itself. Because of its boundedness on the rays $\arg\lambda\not\in\Lambda$,
it follows that $\mathcal{F}(\lambda)\equiv\const$.

Applying Lemma \ref{lm03}, we conclude that $f\equiv0$. We arrived at a contradiction by assuming that there is no completeness with $\rho(\theta)<0$. \qquad$\Box$

\subsection{Step 4. Verification of the sufficient condition}
\begin{Proposition}
\label{prop05}
The function $\rho(\theta)$ increases for $\theta\in [t_0(\alpha),\pi)\cap [t_0(\alpha),\pi\alpha)$, %]$
takes values of different signs in the neighborhoods of the boundaries of this interval.

The only zero of $\rho(\theta)$, $\theta_0(\alpha)=t_0(\alpha)+\Delta t(\alpha)$ lies strictly inside the interval, in particular, $\Delta t(\alpha)>0$, and
is a continuous function for $\alpha\in(0,2)$.

Thus, the sufficient conditions for the completeness of S.E. of the operator  $\mathscr{L}_{c,\alpha}$ are satisfied for all
$\arg c\in [t_0(\alpha),\theta_0(\alpha))$.%]$
\end{Proposition}
{\noindent\bf Proof.} Denote, $\mu_0=e^{it_0(\alpha)}$.

The function $\rho(\theta)$ is continuous at $\theta=t_0(\alpha)$ and $\rho(t_0(\alpha))=-\re S(\zeta_0)<0$, since in this case
$Z_0=\zeta_0$.

Exploring the second boundary value, consider two
cases: $0<\alpha<1$ and $1\le\alpha<2$.

Let $0<\alpha<1$ and $\theta\to\pi\alpha-0$. We substitute the variable in the integral:
$$
S(\zeta_0)=-i\frac{\mu_0^{1/\alpha+1/2}}{c^{1/\alpha}}\int\limits_0^1\sqrt{1-\xi^\alpha}\,d\xi=\frac{i}{c^{1/\alpha}}\int\limits_0^1\sqrt{1-\xi^\alpha}\,d\xi,
$$
Obviously, $\re S(\zeta_0)\to0$, therefore $\rho(\theta)\sim\re S(Z_0)\Bigl|_{\theta=\pi\alpha}\Bigr.>0$.

Let $1\le\alpha<2$ and $\theta\to\pi-0$. Then $\im c\to0$, $\re c\to-1$, $\re c^{1/2}\sim\im c/2$, $Z_0\to+\infty$. In this case $\re S(\zeta_0)$ is bounded.

Let us estimate $\re S(Z_0)$:
\begin{gather*}
\re S(Z_0)>\int\limits_{Z_0/2}^{Z_0}\re\sqrt{c\zeta^\alpha-\mu_0}\,d\zeta=Z_0\int\limits_{1/2}^1\re\sqrt{c\xi^\alpha\frac{\im\mu_0}{\im c}-\mu_0}\,d\xi>\\
>C Z_0\frac{(\im\mu_0)^{1/2}}{(\im c)^{1/2}}\re c^{1/2}\int\limits_{1/2}^1\xi^{\alpha/2}(1+O(\im c))\,d\xi>
\frac{C}{(\im c)^{1/\alpha-1/2}};
\end{gather*}
the first estimate is legal, since $\re S(z)$ is positive and increasing in the interval $[0,Z_0]$. Then we substitute the variable $\zeta=Z_0\xi$;
the explicit expression is used for $Z_0=(\im\mu_0/\im c)^{1/\alpha}$.
The two final inequalities are valid for $\im c\ll1$. The estimation $O(\im c)$ is uniform in $\xi\in[1/2,1]$.

As well as $\alpha<2$, we conclude $\re S(Z_0)\to+\infty$, thus $\rho(\theta)\to+\infty$.

Now we prove the strict increase of $\rho(\theta)$ on the interval $\theta\in (t_0(\alpha),\pi)\cap (t_0(\alpha),\pi\alpha)$. Note,
\begin{gather}
\notag
\rho(\theta)=\re\Bigl\{\int\limits_{\zeta_0}^{Z_0}\sqrt{c\zeta^\alpha-\mu_0}\,d\zeta-\int\limits_0^{\zeta_0}\sqrt{c\zeta^\alpha-\mu_0}\,d\zeta\Bigr\}=\\
\label{eqrhotwoints}
=\re\int\limits_{\zeta_0}^{Z_0}\sqrt{c\zeta^\alpha-\mu_0}\,d\zeta-
\sin\frac{\theta}{\alpha}\int\limits_0^1\sqrt{1-\xi^\alpha}\,d\xi.
\end{gather}

We turn to the first term in \eqref{eqrhotwoints}, denoting
$$
I(\theta)=\re\int\limits_{\zeta_0}^{Z_0}\sqrt{c\zeta^\alpha-\mu_0}\,d\zeta.
$$

For $I(\theta)$ we will carry out the integration along the path $\gamma=\{\zeta=\zeta(t),\ t\in[0,\tau]\}$, $\zeta(t)=((\mu_0-t)/c)^{1/\alpha}$,
$\tau=\im(c/\mu_0)/\im c>0$.

The $(\mu_0-t)/c$ -- image of the interval $t\in[0,\tau]$ is a segment in the IVth quarter (due to the location of endpoints --- Proposition \ref{prop01}).
The IVth quarter also contains the path $\gamma$ --- this is true for endpoints by virtue of the same Proposition \ref{prop01}, and taking into
account the continuity of the $\arg$ function ---
for the whole path $\gamma$.

For the main branch $q^{1/2}(\zeta)=\sqrt{c\zeta^\alpha-\mu_0}=-it^{1/2}(\zeta)$ along $\gamma$. Indeed, according to Proposition \ref{propSonLs}, $\im q^{1/2}(Z_0)<0$,
i.e $q^{1/2}(Z_0)=-i\tau^{1/2}$; for other points of $\gamma$ the sign
is preserved due to the continuity of the square root.

To shorten the notation, we denote the partial derivatives  by indices with respect to the corresponding variables.
For example, $\zeta_t$ is the derivative of the parameterization function for $\gamma$.

We have:
\begin{equation}
\label{eqrhoviazeta}
I(\theta)=\int\limits_0^\tau t^{1/2}\im\zeta_t\,dt,\quad
I_\theta(\theta)=\tau^{1/2}\im\zeta_t(\tau)\tau_\theta+\int\limits_0^\tau t^{1/2}\im\zeta_{\theta t}\,dt.
\end{equation}

Obviously, $c_\theta=ic$, $\tau_\theta=-(\re c/\im c)_\theta\im\mu_0=-(\re c_\theta\im c-\im c_\theta\re c)\im\mu_0/(\im c)^2$, thus
$\tau_\theta=\im\mu_0/(\im c)^2$, as well as $|c|=1$.

Further, $\zeta_t(\tau)=-(1/\alpha)\zeta(\tau)/(\mu_0-\tau)$. As $\zeta(\tau)=Z_0$ and $\mu_0-\tau=cZ_0^\alpha=c\im\mu_0/\im c$, we obtain
$\im\zeta_t(\tau)\tau_\theta=Z_0/\alpha$.
Also note, that $\zeta_\theta=\zeta_c c_\theta=ic\zeta_c=-i\zeta/\alpha$.

Combining these results with \eqref{eqrhoviazeta}, we get:
$$
I_\theta(\theta)=\frac{1}{\alpha}\tau^{1/2}Z_0-\frac{1}{\alpha}\int\limits_0^\tau t^{1/2}\re\zeta_t\,dt.
$$

Applying the second mean value theorem for the integral term, we find some intermediate point $\vartheta\in[0,\tau]$ for which
$$
I_\theta(\theta)=\frac{1}{\alpha}\tau^{1/2}Z_0-\frac{1}{\alpha}\tau^{1/2}(\re\zeta(\tau)-\re\zeta(\vartheta))=\frac{1}{\alpha}\tau^{1/2}\re\zeta(\vartheta)>0,
$$
we have used the facts that all points of $\gamma$ lie strictly in the right half-plane, and $\tau>0$.

Consider the second term in \eqref{eqrhotwoints}:
$$
J(\theta)=\sin\frac{\theta}{\alpha}\int\limits_0^1\sqrt{1-\xi^\alpha}\,d\xi.
$$

Since $t_0(\alpha)<\theta<\pi\alpha$, then $\pi/2<2\pi/(\alpha+2)<\theta/\alpha<\pi$, therefore:
$$
J_\theta(\theta)=\frac{1}{\alpha}\cos\frac{\theta}{\alpha}\int\limits_0^1\sqrt{1-\xi^\alpha}\,d\xi<0,
$$
after all we obtain that $\rho_\theta(\theta)=I_\theta(\theta)-J_\theta(\theta)>0$. The function $\rho(\theta)$ is strictly increasing, has the only zero
$\theta_0=\theta_0(\alpha)=t_0(\alpha)+\Delta t(\alpha)\in (t_0(\alpha),\pi)\cap(t_0(\alpha),\pi\alpha)$. For all  $\theta\in [t_0(\alpha),\theta_0(\alpha))$, %]$
the value $\rho(\theta)<0$.

For any $\alpha_0\in(0,2)$, $\theta_0=\theta_0(\alpha_0)$,
consider $\rho(\theta)=R(\alpha,\theta)$ as a function of two variables: $\alpha$ and $\theta$
in the small neighborhood $\Upsilon$ of
$(\alpha_0,\theta_0)$.

From the representation \eqref{eqrhotwoints} and the first formula of \eqref{eqrhoviazeta} we obtain the existence of
a small neighborhood $\Upsilon=\Upsilon(\alpha_0,\theta_0)$, in which $R(\alpha,\theta)$ is continuous as a function of
two variables.

Since $R(\alpha_0,\theta_0)=0$ and $R(\alpha_0,\theta)$ is strictly monotone as a function of $\theta$ for $(\alpha_0,\theta)\in\Upsilon$,
the Implicit Function Theorem implies the continuity of
$\theta_0(\alpha)$ in some small neighborhood of the point $\alpha_0$. Since $\alpha_0$ is arbitrary, it's valid for all
$\alpha\in(0,2)$.
\qquad$\Box$
\bigskip

We have proved the completeness of the of S.E. of the operator $\mathscr{L}_{c,\alpha}$ in case $t_0(\alpha)\le|\arg c|<\theta_0(\alpha)$. As already noted, the completeness
in case $|\arg c|<t_0(\alpha)$ is a known fact \cite{Savchuk1}. This completes the proof of the Theorem \ref{th01}.

\section{Annex}

We turn to the complex Airy operator to evaluate $\theta_0(1)$ as an example. Let $\alpha=1$; taking into account \eqref{eqrhotwoints}, \eqref{eqrhoviazeta} and the
evaluation
$\zeta_t(t)=-1/c$, we obtain:
$$
\rho(\theta)=\int\limits_0^{\tau(\theta)} t^{1/2}\sin\theta\,dt-\frac{2}{3}\sin\theta=
\frac{2}{3}(\tau^{3/2}(\theta)-1)\sin\theta.
$$

We solve the equation $\rho(\theta_0)=0$ for $\theta_0=\theta_0(1)\in(t_0(1),\pi)=(2\pi/3,\pi)$. The only possible case --- $\tau(\theta_0)=1$.
As well as $\tau(\theta)=\im(c/\mu_0)/\im c=\sin(\theta-2\pi/3)/\sin\theta$, the equation on $\theta_0$ takes the form of:
$$
\sin(\theta_0-2\pi/3)=\sin\theta_0,
$$
and it has the only solution $\theta_0=5\pi/6$, which exactly corresponds to the result of Savchuk and Shkalikov \cite{Savchuk1}.
\bigskip

The author expresses his deep appreciation to Andrei Andreyevich Shkalikov for his attention to the work and valuable advice,
as well as to the team of the scientific seminar ''Operator Models in Mathematical Physics`` for their support.

This research was supported by RFBR grant No 19-01-00240.

\end{document}